\documentclass[12pt]{amsart}
\usepackage{fancyhdr}
\usepackage{amssymb,latexsym,amscd}
\usepackage[english]{babel}
\usepackage{amsthm,amsmath}
\input xy
\xyoption{all}
\newcommand{\cV}{\mathcal{V}}
\newcommand{\cU}{\mathcal{U}}
\newcommand{\cO}{\mathcal{O}}
\newcommand{\cE}{\mathcal{E}}
\newcommand{\cF}{\mathcal{F}}
\newcommand{\cM}{\mathcal{M}}
\newcommand{\cN}{\mathcal{N}}
\newcommand{\cL}{\mathcal{L}}
\newcommand{\cG}{\mathcal{G}}

\newcommand{\cB}{\mathcal{B}}
\newcommand{\cQ}{\mathcal{Q}}
\newcommand{\cK}{\mathcal{K}}
\newcommand{\cD}{\mathcal{D}}
\newcommand{\cH}{\mathcal{H}}

\newcommand{\lra}{\longrightarrow}

\newcommand{\ra}{\rightarrow}
\newcommand{\lms}{\longmapsto}
\newcommand{\ms}{\mapsto}

\newcommand{\gt}{\theta}
\newcommand{\PP}{\mathbb{P}}
\newcommand{\ZZ}{\mathbb{Z}}

\newcommand{\CC}{\mathbb{C}}
\newcommand{\MMM}{\mathfrak{M}}
\newcommand{\fV}{\mathfrak{V}}

\newcommand{\Ext}{\mathrm{Ext}}

\newcommand{\Hom}{\mathrm{Hom}}

\newcommand{\rk}{\mathrm{rk}}
\newcommand{\id}{\mathrm{id}}

\newcommand{\Spec}{\mathrm{Spec}}

\theoremstyle{plain}

\begin{document}
\title[]{On Frobenius-destabilized rank-$2$ vector bundles over curves}
\author{Herbert Lange}
\author{Christian Pauly}
\address{Mathematisches Institut \\ Universit\"at Erlangen-N\"urnberg \\
Bismarckstrasse 1 1/2 \\ D-91054 Erlangen \\ Deutschland}
\email{lange@mi.un-erlangen.de}
\address{D\'epartement de Math\'ematiques \\ Universit\'e de Montpellier II - Case Courrier 051 \\ Place Eug\`ene Bataillon \\ 34095 Montpellier Cedex 5 \\ France}
\email{pauly@math.univ-montp2.fr}
\thanks{}
\begin{abstract}
Let $X$ be a smooth projective curve of genus $g \geq 2$ over an algebraically 
closed field $k$ of characteristic $p > 0$. Let 
$\cM_X$ be the moduli space of semistable rank-2 vector bundles over $X$ with trivial determinant.
The relative Frobenius map $F: X \ra X_1$ induces by pull-back a 
rational map $V: \cM_{X_1} \dashrightarrow \cM_{X}$. In this paper we
show the following results.
\begin{enumerate}
\item For any line bundle $L$ over $X$, the rank-$p$ vector bundle
$F_*L$ is stable.
\item The rational map $V$ has base points, i.e., there exist stable
bundles $E$ over $X_1$ such that $F^* E$ is not semistable.
\item Let $\cB \subset \cM_{X_1}$ denote the scheme-theoretical base locus of $V$. 
If $g=2$, $p>2$ and $X$ ordinary, then $\cB$ is a $0$-dimensional local complete 
intersection of length $\frac{2}{3} p(p^2 -1)$ and the degree of $V$ equals 
$\frac{1}{3}p(p^2+2)$.
\end{enumerate}
\end{abstract}
\maketitle

%%%%%%%%%%%%%%%%%%%%%%%%%%%%%%%%%%%%%%%%%%%%%%%%%%%%%%%%%%

\begin{center}
{\bf Introduction}
\end{center}

\vspace*{0.5cm}

Let $X$ be a smooth projective curve of genus $g \geq 2$ over an algebraically 
closed field $k$ of characteristic $p > 0$. Denote by $F : X \ra X_1$  
the relative $k$-linear Frobenius map. Here $X_1 = X \times_{k,\sigma}k$, 
where $\sigma: \Spec(k) \ra \Spec(k)$ is the Frobenius of $k$ (see e.g. \cite{R} section 4.1).
We denote by $\cM_X$, respectively $\cM_{X_1}$, the moduli
space of semistable rank-2 vector bundles on $X$, respectively $X_1$, with trivial determinant.
The Frobenius $F$ induces by pull-back a rational map (the Verschiebung)
$$ V: \cM_{X_1} \dashrightarrow \cM_{X}, \qquad [E] \ms [F^* E].$$
Here $[E]$ denotes
the S-equivalence class of the semistable bundle $E$.
It is shown \cite{MS} that $V$ is generically \'etale, hence separable
and dominant, if $X$ or equivalently $X_1$ is an ordinary curve. Our first result is

\bigskip

{\bf Theorem 1} 
{\em Over any smooth projective curve $X_1$ of genus $g \geq 2$ there exist stable rank-$2$
vector bundles $E$ with trivial determinant, such that $F^* E$ is not
semistable. In other words, $V$ has base points.}
 
\bigskip

Note that this is a statement for an arbitrary curve of genus $g \geq 2$ over $k$,
since associating $X_1$ to $X$ induces an automorphism of the moduli space of curves
of genus $g$ over $k$. The existence of Frobenius-destabilized bundles was already proved in
\cite{LP2} Theorem A.4 by specializing the so-called Gunning bundle on a Mumford-Tate
curve. The proof given in this paper is much simpler than the previous one.  Given a line
bundle $L$ over $X$, the generalized Nagata-Segre theorem asserts the 
existence of rank-$2$ subbundles $E$ of the rank-$p$ bundle $F_* L$ of 
a certain (maximal) degree. Quite surprisingly, these subbundles $E$ of 
maximal degree turn out to be stable and Frobenius-destabilized. 

\bigskip

In the case $g=2$ the moduli space $\cM_X$ is canonically isomorphic to the 
projective space $\PP_k^3$ and the set of strictly semistable bundles 
can be identified with the Kummer surface $\mathrm{Kum}_X \subset \PP_k^3$ associated to $X$. 
According to \cite{LP2} Proposition A.2 the rational map
$$ V: \PP_k^3 \dashrightarrow \PP_k^3 $$
is given by polynomials of degree $p$, which are explicitly known in the cases $p=2$ \cite{LP1} and 
$p=3$ \cite{LP2}. Let $\cB$ be the scheme-theoretical base locus of $V$, i.e., the
subscheme of $\PP_k^3$ determined by the ideal generated by the $4$ polynomials of degree $p$
defining $V$.  Clearly its
underlying set equals 
$$ \mathrm{supp} \ \cB = \{ E \in  \cM_{X_1} \cong \PP_k^3 \ | \ F^* E \ \ \mbox{is not semistable} \}$$
and $\mathrm{supp} \ \cB  \subset \PP_k^3  \setminus \mathrm{Kum}_{X_1}$.  
Since $V$ has no base points on the ample divisor $\mathrm{Kum}_{X_1}$, we deduce
that $\dim \cB = 0$. Then we show 

\bigskip
{\bf Theorem 2}
{\em Assume $p>2$. Let $X_1$ be an ordinary curve of genus $g=2$. Then the 
$0$-dimensional scheme $\cB$ is a local complete intersection of length}
 $$\frac {2}{3} p(p^2 -1).$$

\bigskip

\noindent
Since $\cB$ is a local complete intersection, the degree of $V$ equals $\deg V = p^3 - l(\cB)$ where $l(\cB)$ denotes the length of $\cB$ (see e.g. \cite{O1} Proposition 2.2). Hence we 
obtain the

\bigskip
{\bf Corollary} {\em Under the assumption of Theorem 2} $$\deg V = \frac{1}{3}p(p^2+2).$$

\bigskip

The underlying idea of the proof of Theorem 2 is rather simple: we observe that
a vector bundle $E \in \mathrm{supp} \ \cB$ corresponds via adjunction to a subbundle
of the rank-$p$ vector bundle $F_*(\gt^{-1})$ for some theta characteristic $\gt$
on $X$ (Proposition 3.1). This is the motivation to introduce Grothendieck's
Quot-Scheme $\cQ$ parametrizing rank-$2$ subbundles of degree $0$ of the
vector bundle $F_*(\gt^{-1})$.  We prove that the two $0$-dimensional schemes
$\cB$ and $\cQ$ decompose as disjoint unions $\coprod \cB_\gt$ and $\coprod
\cQ_\eta$ where $\gt$ and $\eta$ vary over theta characteristics on $X$ and
$p$-torsion points of $JX_1$ respectively and that $\cB_\gt$ and $\cQ_0$ 
are isomorphic, if $X$ is ordinary (Proposition 4.6). In particular 
since $\cQ$ is a local complete intersection, $\cB$ also is.

\bigskip

In order to compute the length of $\cB$ we show that $\cQ$ is isomorphic
to a determinantal scheme $\cD$ defined intrinsically by the $4$-th
Fitting ideal  of some sheaf. The non-existence of a universal family over
the moduli space of rank-$2$ vector bundles of degree $0$ forces us to 
work over a different parameter space constructed via the Hecke 
correspondence and carry out the Chern class computations on this
parameter space.

\bigskip

The underlying set of points of $\cB$ has already been studied in 
the literature. In fact, using the notion of $p$-curvature, S. Mochizuki
\cite{Mo} describes points of $\cB$ as ``dormant atoms'' and 
obtains, by degenerating the genus-$2$ curve $X$ to a singular 
curve, the above mentioned formula for their number (\cite{Mo} Corollary 3.7 
page 267).  Moreover he shows
that for a general curve $X$ the scheme $\cB$ is reduced.
In this context we also mention the recent work of B. Osserman \cite{O2}, \cite{O3},
which explains the relationship of $\mathrm{supp}  \ \cB$ with Mochizuki's theory.

\bigskip
\begin{center}
{\bf Acknowledgments}
\end{center}
\bigskip
We would like to thank Yves Laszlo and Brian Osserman for helpful discussions and
for pointing out several mistakes in a previous version of this paper. We also thank
Adrian Langer for some advice with references.  
  
\bigskip
   
\begin{center}
{\bf \S 1 Stability of the direct image $F_*L$.}
\end{center}
\bigskip

Let $X$ be a smooth projective curve of genus $g \geq 2$ over an algebraically closed field of characteristic $p > 0$ and let
$F: X \ra X_1$ denote the relative Frobenius map.
Let $L$ be a line bundle over $X$ with 
$$\deg L = g-1 + d,$$
for some integer $d$. Applying the Grothendieck-Riemann-Roch theorem
to the morphism $F$, we obtain

\bigskip

{\bf Lemma 1.1} {\it
The slope of the rank-$p$ vector bundle $F_*L$ equals}
$$\mu(F_*L) = g-1 + \frac{d}{p}.$$

The following result will be used in section 3.

\bigskip

{\bf Proposition 1.2}  {\it
If $g \geq 2$, then the vector bundle $F_*L$ is stable for any line bundle $L$ on $X$.}

\begin{proof} 
Suppose that the contrary holds, i.e., $F_*L$ is not stable. 
Consider its Harder-Narasimhan filtration 
$$
0 = E_0 \subset E_1 \subset E_2 \subset \cdots \subset E_l = F_*L,
$$
such that the quotients $E_i/E_{i-1}$ are semistable with $\mu(E_i/E_{i-1}) > \mu(E_{i+1}/E_{i})$ for all $i \in \{ 1,\ldots, l-1 \}$.
If $F_*L$  is not
semistable,  we 
denote $E := E_1$. If $F_*L$ is 
semistable, we denote by $E$ any proper semistable subbundle of the same slope.  Then clearly
\begin{equation}\label{des}
\mu(E) \geq \mu(F_*L).
\end{equation}
In case $r = \rk  \ E > \frac{p-1}{2}$, we observe that the quotient bundle 
$$
Q = \left\{ \begin{array}{ll} F_*L / E_{l-1}& \mbox{if} \; F_*L \; \mbox{is not semistable},\\
                               F_*L / E      & \mbox{if} \; F_*L \; \mbox{is semistable},
                             \end{array} \right.
$$
is also semistable and that its dual
$Q^*$ is a subbundle of $(F_*L)^*$. Moreover, by relative duality $(F_*L)^*= F_*(L^{-1} \otimes
\omega_{X}^{\otimes 1-p})$ and by assumption
$\rk \ Q^* \leq p-r \leq \frac{p-1}{2}$. Hence, replacing
if necessary $E$ and $L$ by $Q^*$ and $L^{-1} \otimes 
\omega_{X}^{\otimes 1-p}$, we may assume that $E$ is semistable
and $r \leq \frac{p-1}{2}$.

\medskip
Now, by \cite{SB} Corollary 2, we have the inequality
\begin{equation} \label{ineq}
\mu_{max}(F^*E) - \mu_{min}(F^*E) \leq (r-1) (2g-2),
\end{equation}
where $\mu_{max}(F^*E)$ (resp. $\mu_{min}(F^*E)$) denotes the slope
of the first (resp. last) graded piece of the Harder-Narasimhan
filtration of $F^*E$. The inclusion $E \subset F_*L$ gives, by adjunction, a nonzero
map $F^*E \ra L$. Hence
$$ \deg L \geq \mu_{min} (F^*E) \geq  \mu_{max} (F^*E) - (r-1)(2g-2) \geq
p\mu(E) - (r-1)(2g-2).$$
Combining this inequality with \eqref{des} and using Lemma 1.1, we obtain
$$ g-1 + \frac{d}{p} = \mu(F_*L) \leq \mu(E) \leq \frac{g-1+d}{p} + \frac{(r-1)(2g-2)}{p},$$
which simplifies to
$$ (g-1) \leq (g-1) \left( \frac{2r-1}{p} \right).$$
This is a contradiction, since we have assumed $r \leq \frac{p-1}{2}$ and therefore
$\frac{2r-1}{p}  < 1$.
\end{proof}

{\bf Remark 1.3} We observe that the vector bundles $F_*L$ are destabilized 
by Frobenius, because of the nonzero canonical map $F^*F_*L \ra L$ and
clearly $\mu(F^*F_*L) > \deg L$. For further properties of the bundles
$F_*L$, see \cite{JRXY} section 5.

\bigskip

{\bf Remark 1.4} In the context of Proposition 1.2 we mention the following open
question: given a finite separable morphism between smooth curves $f: Y \ra X$
and a line bundle $L \in \mathrm{Pic}(Y)$, is the direct image $f_* L$ stable?
For a discussion, see \cite{B}.

\bigskip

\begin{center}
{\bf \S 2 Existence of Frobenius-destabilized bundles.}
\end{center}
\vspace*{0.5cm}
\noindent
Let the notation be as in the previous section.
We recall the generalized Nagata-Segre theorem,
proved by Hirschowitz, which says

\bigskip

{\bf Theorem 2.1}   {\it 
For any vector bundle $G$ of rank $r$ and degree $\delta$ over any smooth
curve $X$ and for any integer $n$, $1 \leq n \leq r-1$, there exists
a rank-$n$ subbundle $E \subset G$, satisfying
\begin{equation}  \label{hb}
\mu(E) \geq \mu(G) - \left(\frac{r-n}{r} \right)(g-1) - \frac{\epsilon}{rn}, 
\end{equation} 
where $\epsilon$ is the unique integer with $0 \leq \epsilon \leq r-1$ and
$\epsilon + n(r-n)(g-1) \equiv n \delta  \ \mathrm{mod} \  r$.}

\bigskip

{\bf Remark 2.2} The previous theorem can be deduced (see \cite{L} Remark 3.14) from the main 
theorem of \cite{H} (for its proof, see http://math.unice.fr/\~{}ah/math/Brill/).

\bigskip

{\bf Proof of Theorem 1.}  We apply Theorem 2.1 to the rank-$p$
vector bundle $F_*L$ on $X_1$ and $n=2$, where $L$ is a line bundle of degree $g-1+d$ on $X$, with
$d \equiv -2g+2 \ \mathrm{mod} \ p$: There exists a rank-$2$ vector bundle
$E \subset F_*L$ such that
\begin{equation}\label{ineq1}
\mu(E)  \geq \mu(F_*L) -  \frac{p-2}{p}  (g-1).
\end{equation}
Note that our assumption on $d$ was made to have $\epsilon = 0$.
 
\medskip
Now we
will check that any $E$ satisfying inequality \eqref{ineq1} is stable with 
$F^* E$ not semistable.

\medskip

(i) $E$ is stable: Let $N$ be a line subbundle of $E$. The inclusion
$N \subset F_*L$ gives, by adjunction, a nonzero map $F^*N \ra L$, which
implies (see also \cite{JRXY} Proposition 3.2(i))
$$ \deg N \leq \mu(F_*L) - \frac{p-1}{p} (g-1).$$
Comparing with \eqref{ineq1} we see that $\deg N < \mu(E)$.
\medskip

(ii) $F^*E$ is not semistable. In fact, we claim that $L$ destabilizes $F^*E$.
For the proof note that Lemma 1.1 implies
\begin{equation} \label{ineq2}
\mu(F_*L) - \frac{p-2}{p} (g-1) = \frac{2g-2+d}{p} > \frac{g-1+d}{p} = \frac{\deg L}{p}
\end{equation}
since $g \geq 2$.  Together with (\ref{ineq1}) this gives
$ \mu(E) > \frac{\deg L}{p}$ and hence 
$$
\mu(F^*E) > \deg L.
$$
This implies the assertion, since by adjunction we obtain a nonzero map $F^*E \ra L$.

\medskip

Replacing $E$ by a subsheaf of suitable degree, we may assume that 
inequality \eqref{ineq1} is an equality. In that case, because of
our assumption on $d$, $\mu(E)$ is an integer, hence $\deg E$ is even.
In order to get trivial determinant, we may tensorize $E$ with a suitable line
bundle. This completes the proof of Theorem 1. \hfill{$\square$}

\bigskip

\begin{center}
{\bf \S 3 Frobenius-destabilized bundles in genus 2.}
\end{center}
\vspace*{0.5cm}

From now on we assume that $X$ is an ordinary  curve of genus $g=2$ and 
the characteristic of $k$ is $p > 2$.
Recall that $\cM_X$  denotes the moduli space of semistable rank-$2$ vector bundles with
trivial determinant over $X$ and $\cB$ the scheme-theoretical base locus of the rational
map
$$ V: \cM_{X_1} \cong \PP_k^3 \dashrightarrow \PP_k^3 \cong \cM_{X},$$
which is given by polynomials of degree $p$.

%  In the sequel we will use the same
% letter to denote line bundles on $X$ and their pull-back to the Frobenius twist $X_1$ under the
% canonical $k$-semilinear isomorphism $\iota : X_1 \rightarrow X$ (see e.g. \cite{R} section 4.1).

\bigskip

First of all we will show that the $0$-dimensional scheme $\cB$ is the disjoint union of subschemes $\cB_\theta$ indexed by theta characteristics of $X$.

\bigskip

{\bf Proposition 3.1} {\it 

\begin{enumerate}

\item[(a)] Let $E$ be a vector bundle on $X_1$ such that $E \in \mathrm{supp} \ \cB$. Then there 
exists a unique theta characteristic 
$\gt$ on $X$, 
such that $E$ is a subbundle of $F_*(\gt^{-1})$.

\item[(b)]
Let $\gt$ be a theta characteristic on $X$. Any rank-$2$ subbundle 
$E \subset F_*(\gt^{-1})$ of degree 0 has the following properties
\begin{itemize}
\item[(i)] $E$ is stable and $F^*E$ is not semistable,
\item[(ii)] $F^*(\det E) = \cO_{X}$,
\item[(iii)] $ \dim \Hom(E,F_*(\gt^{-1}))= 1$ and $ \dim H^1(E^* \otimes F_*(\gt^{-1})) = 5$,
\item[(iv)] $E$ is a rank-$2$ subbundle of maximal degree.
\end{itemize}

\end{enumerate}
}

\bigskip
\noindent
{\it Proof}: (a) By \cite{LS} Corollary 2.6 we know that, for every $E \in \mathrm{supp} \ \cB$
the bundle $F^*E$ is the nonsplit extension of $\gt^{-1}$ by $\gt$, for some theta 
characteristic $\gt$ on $X$ (note that $\Ext^1(\gt^{-1},\gt) \cong k$). By adjunction we get 
a homomorphism $\psi: E \ra F_*(\gt^{-1})$ and we 
have to show that this is of maximal rank. 

Suppose it is not, then there is a line bundle $N$ on the curve $X_1$ such that 
$\psi$ factorizes as   
$E \ra N \ra F_*(\gt^{-1}).$
By stability of $E$ we have $\deg N > 0$. On the other hand, by adjunction, we get a nonzero 
homomorphism $F^*N \ra \gt^{-1}$ implying $p \cdot \deg N \leq -1$, a contradiction.
Hence $\psi: E \ra F_*(\gt^{-1})$ is injective. Moreover $E$ is even a subbundle of 
$F_*(\gt^{-1})$, 
since otherwise there exists a subbundle $E' \subset F_*(\gt^{-1})$ with $\deg E'> 0$ 
and which fits into the exact sequence 
$$ 0 \lra E \lra E' \stackrel{\pi}{\lra} T \lra 0,$$
where $T$ is a torsion sheaf supported on an effective divisor. Varying $\pi$,
we obtain a family of bundles $\ker \pi \subset E'$  of dimension $>0$  and
$\det \ker \pi = \cO_{X_1}$. This would imply (see proof 
of Theorem 1) $\dim \cB > 0$, a contradiction. 

Finally, since $\theta$ is the maximal destabilizing line subbundle of $F^*E$, it is
unique.

\bigskip

(b) We observe that inequality \eqref{ineq1} holds for the pair $E \subset
F_*(\gt^{-1})$. Hence, by the proof of Theorem 1, $E$ is stable and $F^* E$ is not
semistable.

Let $\varphi: F^*E \ra \theta^{-1}$ denote the homomorphism, adjoint to the inclusion $E \subset F_*(\gt^{-1})$.
The homomorphism $\varphi$ is surjective, since otherwise $F^*E$ would contain a line subbundle of degree $>1$, contradicting 
\cite{LS}, Satz 2.4. Hence we get an exact sequence
\begin{equation} \label{es}
0 \ra \ker \varphi \ra F^*E \ra \gt^{-1} \ra 0.
\end{equation}

On the other hand, let $N$ denote a line bundle on $X_1$ such that $E \otimes N$ has trivial determinant, i.e. $N^{-2} = \det E$.
Applying \cite{LS} Corollary 2.6 to the bundle $F^*(E \otimes N)$ we get an exact sequence
$$
0 \ra \tilde{\gt} \otimes N^{-p} \ra F^*E \ra \tilde{\gt}^{-1} \otimes N^{-p} \ra 0,
$$
for some theta characteristic  $\tilde{\gt}$.
By uniqueness of the destabilizing subbundle of maximal degree of $F^*E$, this
exact sequence
must coincide with (\ref{es}) up to a nonzero constant. This implies that
$N^p \otimes \tilde{\gt} = \gt$, hence $N^{2p} = \cO_X$. So we obtain that 
$\cO_{X} = \det (F^*E) = F^*(\det E)$ proving (ii).

\bigskip

By adjunction we have
the equality $\dim \Hom(E, F_*(\gt^{-1})) = \dim \Hom(F^* E, \gt^{-1}) = 1$. 
Moreover by Riemann-Roch we obtain $\dim H^1(E^* \otimes F_*(\gt^{-1})) = 5$. This
proves (iii). 

\bigskip

Finally, suppose that there exists a rank-$2$ subbundle $E' \subset F_*(\gt^{-1})$
with $\deg E' \geq 1$. Then we can consider  the kernel $E = \ker \pi$ of a
surjective morphism $\pi : E' \ra T$ onto a torsion sheaf with length equal to
$\deg E'$. By varying $\pi$ and after tensoring $\ker \pi$ with a suitable
line bundle of degree $0$, we construct a family of dimension $>0$ of stable
rank-$2$ vector bundles with trivial determinant which are
Frobenius-destabilized, contradicting $\dim \cB =0$. This proves (iv).
\qed \\

It follows from Proposition 3.1 (a) that the scheme $\cB$ decomposes as a disjoint 
union
$$\cB = \coprod_{\gt} \cB_\gt,$$
where $\theta$ varies over the set of all theta characteristics of $X$ and
$$ \mathrm{supp} \ \cB_\gt = \{ E \in \mathrm{supp} \ \cB  \,\,| \,\, E \subset F_*(\gt^{-1}) \}.$$

\bigskip

Tensor product with a $2$-torsion point $\alpha \in JX_1[2] \cong JX[2]$ induces an isomorphism
of $\cB_{\gt}$ with $\cB_{\gt \otimes \alpha}$ for every theta characteristic $\gt$. 
We denote by $l(\cB)$ and $l(\cB_\gt)$ the length of the schemes $\cB$ and $\cB_\gt$. From the 
preceding we deduce the relations 
\begin{equation} \label{ls}
l(\cB) = 16 \cdot l(\cB_\gt) \qquad \text{for every theta characteristic} \ \ \gt.
\end{equation}

\newpage

\begin{center}
{\bf \S 4 Grothendieck's Quot-Scheme.}
\end{center}
\vspace*{0.5cm}

Let  $\gt$ be a theta characteristic on $X$. We consider the functor $\underline{\cQ}$
from the opposite category of $k$-schemes to the category of sets defined by

\begin{align*}
\underline{\cQ}(S)= \{  & \sigma : \pi_{X_1}^*(F_*(\gt^{-1})) \rightarrow \cG \rightarrow 0 \ |  \ 
\cG \ \text{coherent over} \  X_1 \times S, \ \text{flat over} \  S, \\
& \deg \cG_{|X_1 \times \{ s \} } = \mathrm{rk} \  \cG_{|X_1 \times \{ s \}} = 
p-2 , \ \forall s \in S \} / \cong
\end{align*}

where $\pi_{X_1} : X_1 \times S \rightarrow X_1$ denotes the natural projection and
$\sigma \cong \sigma'$ for quotients $\sigma$ and $\sigma'$ if and only if there
exists an isomorphism $\lambda: \cG \rightarrow \cG'$ such that $\sigma' = \lambda \circ
\sigma$.

\bigskip

Grothendieck showed in \cite{G} (see also \cite{HL} section 2.2)
that the functor $\underline{\cQ}$ is representable by a $k$-scheme, which
we denote by $\cQ$. A $k$-point of $\cQ$ corresponds to a quotient
$\sigma : F_*(\gt^{-1}) \rightarrow G$, or equivalently to a rank-$2$ 
subsheaf $E = \ker \sigma \subset F_*(\gt^{-1})$ of degree $0$ on $X_1$. By the
same argument as in Proposition 3.1 (a), any subsheaf $E$ of degree $0$ is
a subbundle of $F_*(\gt^{-1})$. Since by Proposition 3.1 (b) (iv) the bundle $E$ has maximal
degree as a subbundle of $F_*(\gt^{-1})$, any sheaf $\cG \in \underline{\cQ}(S)$ is 
locally free (see \cite{MuSa} or \cite{L} Lemma 3.8).

\bigskip

Hence taking the kernel of $\sigma$ induces a bijection of $\underline{\cQ}(S)$ with the
following set, which we also denote by $\underline{\cQ}(S)$

\begin{align*}
\underline{\cQ}(S)= \{  & \cE  \hookrightarrow \pi_{X_1}^*(F_*(\gt^{-1})) \ |  \ 
\cE \ \text{locally free sheaf over} \  X_1 \times S \ \text{of rank} \ 2,  \\
&  \pi_{X_1}^*(F_*(\gt^{-1}))/ \cE \ \text{locally free} \ , \ 
\deg \cE_{|X_1 \times \{ s \} } = 0  \
\forall s \in S \} / \cong
\end{align*}
\bigskip

By Proposition 3.1 (b) the scheme $\cQ$ decomposes as a disjoint union 
$$\cQ = \coprod_\eta \cQ_\eta,$$
where $\eta$ varies over the $p$-torsion points $\eta \in JX_1[p]_{red} = \ker (
V: JX_1 \rightarrow JX)$. We also denote by $V$ the Verschiebung of $JX_1$, i.e. 
$V(L) = F^* L$, for $L \in JX_1$. The set-theoretical support of $\cQ_\eta$ 
equals
$$ \mathrm{supp} \ \cQ_\eta = \{ E \in \mathrm{supp} \ \cQ  \,\,| \,\, 
\det E = \eta \}.$$
Because of the projection formula, tensor product with a $p$-torsion point 
$\beta  \in JX_1[p]_{red}$ induces an isomorphism of $\cQ_\eta$ with
$\cQ_{\eta \otimes \beta}$. So the scheme $\cQ$ is a principal homogeneous
space for the group $JX_1[p]_{red}$ and we have the relation
\begin{equation} \label{lq}
l(\cQ) = p^2 \cdot l(\cQ_0),
\end{equation}
since $X_1$ is assumed to be ordinary. Moreover, by Proposition 3.1 we have
the set-theoretical equality
$$ \mathrm{supp} \ \cQ_0 = \mathrm{supp} \ \cB_\theta.$$

\bigskip

{\bf Proposition 4.1} {\it  

\begin{enumerate}

\item[(a)]  $\dim \cQ = 0$.

\item[(b)]  The scheme $\cQ$ is a local complete intersection at any $k$-point
$e = (E \subset F_*(\gt^{-1})) \in \cQ$.

\end{enumerate}
}

\bigskip
\noindent
{\it Proof}: Assertion (a) follows from the preceding remarks and $\dim \cB = 0$. By
\cite{HL} Proposition 2.2.8 assertion (b) follows from the equality $\dim_{[E]} \cQ
= 0 = \chi(\underline{\mathrm{Hom}}(E,G))$, where $E = \ker ( \sigma : F_*(\gt^{-1}) 
\rightarrow G)$ and $\underline{\mathrm{Hom}}$ denotes the sheaf of homomorphisms.
\qed \\

\bigskip

Let $\cN_{X_1}$ denote the moduli
space of semistable rank-$2$ vector bundles of degree $0$ over $X_1$. 
We denote by $\cN_{X_1}^s$ and $\cM_{X_1}^s$ the open subschemes of 
$\cN_{X_1}$ and $\cM_{X_1}$ corresponding to stable vector bundles.
Recall (see  \cite{La1} Theorem 4.1) that $\cN_{X_1}^s$ and
$\cM_{X_1}^s$  universally corepresent the functors (see e.g. \cite{HL} Definition 
2.2.1) from the opposite
category of $k$-schemes of finite type to the category of sets defined by

\begin{align*}
\underline{\cN}^s_{X_1}(S)= \{ &  \cE
\ \text{locally free sheaf over} \  X_1 \times S \ \text{of rank} \ 2 \ | \ 
\cE_{| X_1 \times \{s\} } \ \text{stable} , \\ & \deg \cE_{| X_1 \times \{s\} } = 0, \ \forall
s \in S \} / \sim,
\end{align*}
\begin{align*}
\underline{\cM}^s_{X_1}(S)= \{ & \cE
\ \text{locally free sheaf over} \  X_1 \times S \ \text{of rank} \ 2 \ | \ 
\cE_{| X_1 \times \{s\} } \ \text{stable} \ \forall
s \in S, \\  &  \det \cE = \pi^*_S M \
 \text{for some line bundle} \ M \ \text{on} \  S \} / \sim,
\end{align*}
where $\pi_S : X_1 \times S \rightarrow S$ denotes the natural projection and $\cE' \sim \cE$ 
if and only if there exists a line bundle $L$ on $S$ such that $\cE' \cong \cE \otimes
\pi_S^* L$. We denote by $\langle \cE \rangle$ the equivalence class of the vector
bundle $\cE$ for the relation $\sim$.

\bigskip

Consider the determinant morphism
$$ \det : \cN_{X_1}  \rightarrow JX_1, \qquad [E] \mapsto \det E,$$
and denote by $\det^{-1}(0)$ the scheme-theoretical fibre over the trivial line
bundle on $X_1$. Since $\cN_{X_1}^s$  universally corepresents the functor 
$\underline{\cN}^s_{X_1}$ , we have an isomorphism 
$$\cM_{X_1}^s \cong \cN_{X_1}^s \cap \mathrm{det}^{-1}(0).$$

\bigskip

{\bf Remark 4.2} If $p>0$, it is not known whether the canonical morphism
$\cM_{X_1} \ra \det^{-1}(0)$ is an isomorphism (see e.g. \cite{La2} section 3).

\bigskip

In the sequel we need the following relative version of Proposition 3.1 (b)(ii). By
a $k$-scheme we always mean a $k$-scheme of finite type.
\bigskip

{\bf Lemma 4.3}  \  {\it Let $S$ be a connected $k$-scheme and let $\cE$ be a locally free 
sheaf of rank-$2$ over
$X_1 \times S$  such that $\deg \cE_{|X_1 \times \{ s\} } = 0$ for all 
points $s$ of $S$. Suppose that $\mathrm{Hom}(\cE,\pi^*_{X_1}(F_*(\gt^{-1})) \not= 0$.
Then we have the exact sequence
$$ 0 \lra \pi_X^*(\gt) \lra ( F \times \mathrm{id}_S)^* \cE \lra \pi_X^*(\gt^{-1}) \lra 0. $$
In particular}
$$ (F \times \mathrm{id}_S)^*(\det \cE) = \cO_{X_1 \times S}.$$

\bigskip

{\it Proof:} First we note that by flat base change for $\pi_{X_1} : X_1 \times 
S \rightarrow X_1$, we have an isomorphism $\pi_{X_1}^* (F_*(\gt^{-1})) \cong
( F \times \mathrm{id}_S)_* (\pi_X^*(\gt^{-1}))$. Hence the nonzero morphism 
$\cE \rightarrow \pi_{X_1}^*(F_*(\gt^{-1}))$ gives via adjunction 
a nonzero morphism
$$ \varphi: ( F \times \mathrm{id}_S)^* \cE \lra \pi_X^*(\gt^{-1}).$$
We know by the proof of Proposition 3.1 (b) that the fibre $\varphi_{(x,s)}$ over any
closed point $(x,s) \in X \times S$  is a surjective $k$-linear map. Hence 
$\varphi$ is surjective by
Nakayama and we have an exact sequence
$$ 0 \lra \cL \lra ( F \times \mathrm{id}_S)^* \cE \lra \pi_X^*(\gt^{-1}) \lra 0,$$
with $\cL$ locally free sheaf of rank $1$.
By \cite{K} section 5, the rank-$2$ vector bundle $(F \times \mathrm{id}_S)^* \cE$ is
equipped with a canonical connection 
$$ \nabla : (F \times \mathrm{id}_S)^* \cE \lra (F \times \mathrm{id}_S)^* \cE 
\otimes \Omega^1_{X \times S/S}.$$
We note that $\Omega^1_{X \times S/S} = \pi^*_X(\omega_X)$, where $\omega_X$ denotes the
canonical line bundle of $X$. The first fundamental form of the connection $\nabla$ is an
$\cO_{X \times S}$-linear homomorphism
$$ \psi_\nabla : \cL \lra \pi^*_X(\gt^{-1}) \otimes \pi^*_X(\omega_X) = \pi^*_X (\gt).$$
The restriction of $\psi_\nabla$ to the curve $X \times \{ s \} \subset X \times S$ for any 
closed point $s \in S$ is an isomorphism (see e.g. proof of \cite{LS} Corollary 2.6). Hence 
the fibre of $\psi_\nabla$ is a $k$-linear isomorphism over any closed point $(x,s) 
\in X \times S$. We conclude that $\psi_\nabla$ is an isomorphism, by
Nakayama's lemma and because $\cL$ is a locally free sheaf of rank $1$.

\bigskip

We obtain the second assertion of the lemma, since
$$(F \times \mathrm{id}_S)^*(\det \cE) = \det (F \times \mathrm{id}_S)^* \cE = 
\cL \otimes \pi_X^*(\gt^{-1}) = \cO_{X_1 \times S}. $$
 \qed \\

\bigskip

{\bf Proposition 4.4}  {\it We assume $X$ ordinary.

\begin{enumerate}
\item[(a)] The forgetful morphism 
$$ i: \cQ \hookrightarrow \cN_{X_1}^s, \qquad e = (E \subset F_*(\gt^{-1})) \mapsto
E $$
is a closed embedding.

\item[(b)] The restriction $i_0$ of $i$ to the subscheme $\cQ_0 \subset \cQ$ factors through
$\cM_{X_1}^s$, i.e. there is a closed embedding
$$i_0 : \cQ_0 \hookrightarrow \cM_{X_1}^s.$$

\end{enumerate}
}

\bigskip
\noindent
{\it Proof:} (a) Let $e= (E \subset F_*(\gt^{-1}))$ be a $k$-point of $\cQ$.  To show that 
$i$ is a closed embedding at $e \in \cQ$, it is enough to show that the differential
$(di)_e : T_e \cQ \rightarrow T_{[E]} \cN_{X_1}$ is injective. The Zariski tangent spaces
identify with $\mathrm{Hom}(E,G)$ and $\mathrm{Ext}^1(E,E)$ respectively (see e.g.
\cite{HL} Proposition 2.2.7 and Corollary 4.5.2). Moreover, if we apply the functor
$\mathrm{Hom}(E, \cdot)$ to the exact sequence associated with $e \in \cQ$
$$ 0 \lra E \lra F_*(\gt^{-1}) \lra G \lra 0,$$
the coboundary map $\delta$ of the long exact sequence
$$ 0 \lra \mathrm{Hom}(E,E) \lra \mathrm{Hom}(E, F_*(\gt^{-1})) \lra \mathrm{Hom}(E,G)
\stackrel{\delta}{\lra} \mathrm{Ext}^1(E,E) \lra \cdots$$
identifies with the differential $(di)_e$. Now since the bundle $E$ is stable,
we have $k \cong \mathrm{Hom}(E,E)$. By Proposition 3.1 (b) we obtain that
the map $\mathrm{Hom}(E,E) \rightarrow \mathrm{Hom}(E, F_*(\gt^{-1}))$ is an
isomorphism. Thus $(di)_e$ is injective.
\bigskip

(b) We consider the composite map
$$\alpha : \cQ \stackrel{i}{\lra} \cN_{X_1}^s \stackrel{\det}{\lra} JX_1 \stackrel{V}{\lra}
JX,$$
where the last map is the isogeny given by the Verschiebung on $JX_1$, i.e. $V(L) = 
F^*L$ for $L \in JX_1$. The morphism $\alpha$ is induced by the natural transformation
of functors $\underline{\alpha} : \underline{\cQ} \Rightarrow \underline{JX}$, defined
by 
$$ \underline{\cQ}(S) \lra \underline{JX}(S), \qquad (\cE \hookrightarrow 
\pi^*_X(F_*(\gt^{-1}))) \mapsto (F \times \mathrm{id}_S)^*(\det \cE).$$
Using Lemma 4.3 this immediately implies that $\alpha$ factors through the 
inclusion of the reduced point $\{ \cO_X \} \hookrightarrow JX$. Hence
the image of $\cQ$ under the composite morphism $\det \circ i$ is
contained in the kernel of the isogeny $V$, which is the reduced scheme
$JX_1[p]_{red}$, since we have assumed $X$ ordinary.
Taking connected components we see that the image of $\cQ_0$ under $\det \circ i$
is the reduced point $\{ \cO_{X_1} \} \hookrightarrow JX_1$, which implies that
$i_0(\cQ_0)$ is contained in $\cN_{X_1}^s \cap \det^{-1}(0) \cong \cM_{X_1}^s$.
\qed \\

\bigskip

In order to compare the two schemes $\cB_\gt$ and $\cQ_0$ we need the
following lemma.

\bigskip

{\bf Lemma 4.5} 
{\it 
\begin{enumerate}
\item  The closed subscheme $\cB \subset \cM_{X_1}^s$ corepresents the functor
$\underline{\cB}$ which associates to a $k$-scheme $S$ the set

\begin{align*}
\underline{\cB}(S) = \{ & \cE  \ \text{locally free sheaf over} \ X_1 \times S \ \text{of
rank} \ 2 \ | \ \cE_{| X_1 \times \{s\} } \ \text{stable} \ \forall
s \in S, \\
& 0 \ra \cL \ra (F \times \mathrm{id}_S)^* \cE \ra \cM \ra 0,  \ \text{for some
locally free sheaves} \ \cL, \cM \\ 
& \text{over} \ X \times S \  \text{of rank} \ 1, 
 \deg \cL_{| X \times \{s\} }  = - \deg \cM_{| X \times \{s\} } = 1 \ \forall s \in S,  \\
& \det \cE = \pi^*_S M \
 \text{for some line bundle} \ M \ \text{on} \  S \} / \sim. \\
\end{align*}

\item The closed subscheme $\cB_\gt \subset \cM_{X_1}^s$ corepresents the subfunctor 
$\underline{\cB}_\gt$ of $\underline{\cB}$ defined by $ \langle \cE \rangle 
\in \underline{\cB}_\gt(S)$
if and only if the set-theoretical image of the classifying morphism of $\cL$
$$ \Phi_{\cL}: S \lra \mathrm{Pic}^1(X), \qquad s \lms \cL_{|X \times \{ s \} },$$
is the point $\theta \in  \mathrm{Pic}^1(X)$.

\end{enumerate}
}

\bigskip

{\it Proof:}  We denote by $\MMM_{X_1}$ the algebraic stack parametrizing
rank-$2$ vector bundles  with trivial determinant over $X_1$. Let 
$\MMM^{ss}_{X_1}$ denote the open substack  of $\MMM_{X_1}$ 
parametrizing semistable bundles and $\MMM_{X_1}^{unst}$ the
closed substack of $\MMM_{X_1}$ parametrizing non-semistable 
bundles. We will use the following facts about the stack
$\MMM_{X_1}$.
\begin{itemize}
\item The pull-back of $\cO_{\PP^3}(1)$ by the natural map 
$\MMM_{X_1}^{ss} \ra \cM_{X_1} \cong \PP^3$ extends to a 
line bundle, which we denote by $\cO(1)$, over the moduli stack
$\MMM_{X_1}$ and $\mathrm{Pic}(\MMM_{X_1}) = \ZZ \cdot \cO(1)$.
Moreover there are natural isomorphisms
$H^0(\MMM_{X_1}, \cO(n)) \cong H^0(\cM_{\PP^3}, \cO_{\PP^3}(n))$
for any positive integer $n$ (see \cite{BL} Propositions 8.3 and 8.4).

\item The closed subscheme $\MMM_{X_1}^{unst}$ is the base locus
of the linear system $|\cO(1)|$ over the stack $\MMM_{X_1}$.  
This is seen as follows: we deduce from \cite{S} Theorem 6.2 that
$\MMM_{X_1}^{unst}$ is the base locus of the linear system 
$|\cO(n)|$  for some integer $n$. Since $|\cO(n)|$ is generated by symmetric products
of $n$ sections in $|\cO(1)|$, we obtain that $\MMM_{X_1}^{unst}$ is the
base locus of $|\cO(1)|$.

\end{itemize}

We need to compute the fibre product functor $\underline{\cB} = \cB 
\times_{\cM_{X_1}^s} \underline{\cM}_{X_1}^s$. Let $\fV : \MMM_{X_1} \ra
\MMM_X$ denote the morphism of stacks induced by pull-back under the
Frobenius map $F: X \ra X_1$. Let $S$ be a $k$-scheme and consider 
a vector bundle $\cE \in \MMM_{X_1}^s(S)$. Since the subscheme $\cB$ is
defined as base locus of the linear system $V^{*} |\cO_{\PP^3}(1)|$, we obtain
that $\langle \cE \rangle  \in \underline{\cB}(S)$ if and only if $\cE$ lies
in the base locus of $\fV^* |\cO(1)|$ --- here we use the
isomorphism $|\cO_{\PP^3}(1)| \cong |\cO(1)|$ ---, or 
equivalently $\fV(\cE) = (F \times \mathrm{id}_S)^* \cE$ lies in the
base locus of $|\cO(1)|$, which is the closed substack $\MMM_{X_1}^{unst}$.

\bigskip

By \cite{Sh} section 5 the substack  $\MMM_{X_1}^{unst,1}$ of $\MMM_{X_1}^{unst}$ parametrizing
non-semistable vector bundles having a maximal destabilizing line subbundle of degree $1$
is an open substack of  $\MMM_{X_1}^{unst}$. By \cite{LS} Corollary 2.6 the 
vector bundle $\fV(\cE)$ lies in $\MMM_{X_1}^{unst,1}(S)$. We then consider the
universal exact sequence defined by the Harder-Narasimhan filtration over
$\MMM_{X_1}^{unst,1}$:
$$ 0 \ra \cL \ra (F \times \mathrm{id}_S)^* \cE \ra \cM \ra 0,$$
with $\cL$ and $\cM$ locally free sheaves over $X \times S$ such that
$\deg \cL_{| X \times \{s\} }  = - \deg \cM_{| X \times \{s\} } = 1$ for any $s \in S$.
This proves (1).

\bigskip

As for (2), we add the condition that the family $\cE$ is Frobenius-destabilized by the
theta-characteristic $\gt$. 
\qed 

\bigskip

{\bf Proposition 4.6} {\it There is a scheme-theoretical equality
$$ \cB_\gt = \cQ_0$$
as closed subschemes of $\cM_{X_1}$.}

\bigskip

{\it Proof:}  Since $\cB_\gt$ and $\cQ_0$ corepresent the two functors $\underline{\cB}_\gt$
and $\underline{\cQ}_0$ it will be enough to show that there is a 
canonical bijection between the set $\underline{\cB}_\gt(S)$ 
and $\underline{\cQ}_0(S)$ for any $k$-scheme $S$. We recall that
\begin{align*}
\underline{\cQ}_0(S)= \{  & \cE  \hookrightarrow \pi_{X_1}^*(F_*(\gt^{-1})) \ |  \ 
\cE \ \text{locally free sheaf over} \  X_1 \times S \ \text{of rank} \ 2,  \\
& \pi_X^*(F_*(\gt^{-1}))/ \cE \ \text{locally free} , \ \det \cE \cong \cO_{X_1 \times S} \}
 / \cong
\end{align*}
Note that the property $\det \cE \cong \cO_{X_1 \times S}$ is implied as follows: by
Proposition 4.4 (b) we have $\det \cE \cong \pi_S^* L$ for some line bundle $L$ over $S$ and
by Lemma 4.3 we conclude that $L = \cO_S$.

\bigskip

First we show that the natural map $\underline{\cQ}_0(S) \lra \underline{\cM}_{X_1}^s(S)$
is injective.  Suppose that there exist $\cE,\cE' \in 
\underline{\cQ}_0(S)$ such that $\langle \cE \rangle = \langle \cE' \rangle$, i.e. 
$\cE' \cong \cE \otimes \pi_S^*(L)$ for some line
bundle $L$ on $S$. Then by Lemma 4.3 we have two inclusions
$$ i: \pi_X^*(\gt) \lra (F \times \mathrm{id}_S)^* \cE, \qquad 
i': \pi_X^*(\gt) \otimes \pi_S^*(L^{-1})  \lra (F \times \mathrm{id}_S)^* \cE.$$
Composing with the projection $\sigma : (F \times \mathrm{id}_S)^* \cE \ra \pi_X^*(\gt^{-1})$
we see that the composite map $\sigma \circ i'$ is identically zero. Hence
$\pi_S^*(L) = \cO_{X_1 \times S}$.

\bigskip

Therefore the two sets $\underline{\cQ}_0(S)$ and $\underline{\cB}_\gt(S)$ are 
naturally subsets of  $\underline{\cM}_{X_1}^s(S)$.

\bigskip

We now show that $\underline{\cQ}_0(S) \subset \underline{\cB}_\gt(S)$. Consider
$\cE \in \underline{\cQ}_0(S)$. By Proposition 3.1 (b) the bundle $\cE_{|X_1 \times
\{ s \}}$ is stable for all $s \in S$. By Lemma 4.3 we can take $\cL = \pi_X^*(\gt)$ and
$\cM = \pi_X^*(\gt^{-1})$, so that $ \langle \cE \rangle \in \underline{\cB}_\gt(S)$.

\bigskip
Hence it remains to show that $\underline{\cB}_\gt(S) \subset \underline{\cQ}_0(S)$. Consider a sheaf $\cE$ with 
$\langle \cE \rangle  \in \underline{\cB}_\gt(S)$ --- see Lemma 4.5 (2). 
As in the proof of Lemma 4.3 we consider 
the canonical connection $\nabla$ on $(F \times \mathrm{id}_S)^* \cE$. Its
first fundamental form is an $\cO_{X \times S}$-linear homomorphism
$$\psi_\nabla : \cL \lra \cM \otimes \pi_X^*(\omega_X),$$
which is surjective on closed points $(x,s) \in X \times S$. Hence we can
conclude that $\psi_\nabla$ is an isomorphism. Moreover taking the
determinant, we obtain
$$ \cL \otimes \cM = \det (F \times \mathrm{id}_S)^* \cE = \pi_S^* M .$$
Combining both isomorphisms we deduce that
$$ \cL \otimes \cL = \pi_X^*(\omega_X) \otimes \pi_S^* M .$$
Hence its classifying morphism 
$\Phi_{\cL \otimes \cL} : S \lra \mathrm{Pic}^2(X)$
factorizes through the inclusion of the reduced point $\{ \omega_X \} \hookrightarrow
\mathrm{Pic}^2(X)$. Moreover the composite map of $\Phi_{\cL}$ with the
duplication map $[2]$
$$\Phi_{\cL \otimes \cL} : S \stackrel{\Phi_\cL}{\lra} \mathrm{Pic}^1(X) 
\stackrel{[2]}{\lra} \mathrm{Pic}^2(X)$$
coincides with $\Phi_{\cL \otimes \cL}$. We deduce that $\Phi_\cL$ factorizes through
the inclusion of the reduced point $\{ \theta \} \hookrightarrow \mathrm{Pic}^1(X)$.
Note that the fibre $[2]^{-1} (\omega_X)$ is reduced, since $p>2$. Since $\mathrm{Pic}^1(X)$ is a fine moduli space, 
there exists a line bundle $N$ over $S$ such that
$$\cL = \pi_X^*(\gt) \otimes \pi_S^*(N).$$
We introduce the vector bundle $\cE_0  = \cE \otimes \pi_S^*(N^{-1})$. Then 
$\langle \cE_0 \rangle = \langle \cE \rangle$ and we have an exact sequence
$$ 0 \lra \pi_X^*(\gt) \lra (F \times \mathrm{id}_S)^* \cE_0 \stackrel{\sigma}{\lra}
\pi_X^*(\gt^{-1}) \lra 0,$$
since $\pi_S^* M = \pi_S^* N^2$.
By adjunction the morphism $\sigma$ gives a nonzero morphism
$$ j: \cE_0 \lra (F \times \mathrm{id}_S)_* (\pi^*_X(\gt^{-1})) \cong
\pi^*_{X_1} (F_*(\gt^{-1})).$$
We now show that $j$ is injective. Suppose it is not. Then there
exists a subsheaf $\tilde{\cE}_0 \subset \pi^*_{X_1} (F_*(\gt^{-1}))$ and
a surjective map $\tau : \cE_0 \ra \tilde{\cE}_0$. Let $\cK$ denote the
kernel of $\tau$. Again by adjunction we obtain a map $\alpha : 
(F \times \mathrm{id}_S)^* \tilde{\cE}_0 \ra \pi_X^*(\gt^{-1})$ such that the
composite map
$$ \sigma : (F \times \mathrm{id}_S)^* \cE_0 \stackrel{\tau^*}{\lra}
 (F \times \mathrm{id}_S)^* \tilde{\cE}_0  \stackrel{\alpha}{\lra}  \pi_X^*(\gt^{-1}) $$
coincides with $\sigma$.  Here $\tau^*$ denotes the map $(F \times \mathrm{id}_S)^* \tau$.
Since $\sigma$ is surjective, $\alpha$ is also surjective.
We denote by $\cM$ the kernel of $\alpha$. The induced map $\overline{\tau} :
\pi_{X}^*(\gt) = \ker \sigma \ra \cM$ is surjective, because $\tau^*$ is
surjective. Moreover the first 
fundamental form of the canonical connection $\tilde{\nabla}$ on 
$(F \times \mathrm{id}_S)^* \tilde{\cE}_0$ induces an $\cO_{X \times S}$-linear
homomorphism $\psi_{\tilde{\nabla}} : \cM \ra \pi_X^*(\gt)$ and
the composite map
$$ \psi_\nabla : \pi_X^*(\gt) \stackrel{\overline{\tau}}{\lra} \cM 
\stackrel{\psi_{\tilde{\nabla}}}{\lra} \pi_X^*(\gt)$$
coincides with the first fundamental form of $\nabla$ of 
$(F \times \mathrm{id}_S)^* \cE_0$, which is an isomorphism. Therefore
$\overline{\tau}$ is an isomorphism too. So $\tau^*$ is an isomorphism and
$(F \times \mathrm{id}_S)^* \cK = 0$. We deduce that $\cK = 0$. 

\bigskip

In order to show that $\cE_0 \in \underline{\cQ}_0(S)$, it remains to verify that 
the quotient sheaf $\pi_{X_1}(F_*(\gt^{-1}))/ \cE_0$ is flat over $S$. We recall that
flatness implies locally freeness because of maximality of degree. But 
flatness follows from \cite{HL} Lemma 2.1.4, since the restriction of $j$ to
$X_1 \times \{ s \}$ is injective for any closed $s \in S$ by Proposition 3.1 (a).  \qed \\

\bigskip
Combining this proposition with relations \eqref{ls} and \eqref{lq}, we obtain

\bigskip

{\bf Corollary 4.7} {\it We have }
$$l(\cB) = \frac{16}{p^2} \cdot l(\cQ).$$

\bigskip

\begin{center}
{\bf \S 5 Determinantal subschemes.}
\end{center}
\vspace*{0.5cm}

In this section we introduce a determinantal subscheme $\cD \subset \cN_{X_1}$,
whose length will be computed in the next section. We also show that $\cD$ is
isomorphic to Grothendieck's Quot-scheme $\cQ$.
We first define a determinantal subscheme $\tilde{\cD}$ of a variety $JX_1 \times Z$
covering $\cN_{X_1}$ and then we show that $\tilde{\cD}$ is a $\PP^1$-fibration  over
an \'etale cover of $\cD \subset \cN_{X_1}$.

\bigskip

Since there does not exist a universal bundle over $X_1 \times \cM_{X_1}$, following
an idea of Mukai \cite{Mu}, we consider the 
moduli space $\cM_{X_1}(x)$ of stable rank-$2$ vector bundles on $X_1$
with determinant $\cO_{X_1}(x)$ for a fixed point $x \in X_1$. 
According to \cite{N1} the variety $\cM_{X_1}(x)$ 
is a smooth intersection of two quadrics in $\PP^5$. Let $\cU$ denote a universal bundle
on $X_1 \times \cM_{X_1}(x)$ and denote 
$$\cU_x := \cU_{|\{x\} \times \cM_{X_1}(x)}$$ 
considered as a rank-$2$ vector bundle on $\cM_{X_1}(x)$. Then the 
projectivized bundle 
$$
Z := \PP(\cU_x)
$$
is a $\PP^1$-bundle over $\cM_{X_1}(x)$. 
The variety $Z$
parametrizes pairs $(F_z,l_z)$ consisting of a stable vector bundle $F_z \in \cM_{X_1}(x)$ and a 
linear form $l_z :F_z(x) \ra k_x$ on the fibre of $F_z$ over $x$. Thus to any 
$z \in Z$ one can associate an exact sequence
$$
0 \ra E_z \ra F_z \ra k_x \ra 0
$$
uniquely determined up to a multiplicative constant. Clearly
$E_z$ is semistable, since $F_z$ is stable, and $\det E_z = \cO_{X_1}$. 
Hence we get a diagram (the so-called
Hecke correspondence)
$$
\xymatrix{
Z \ar[d]_{\pi} \ar[r]^(.3){\varphi} & \cM_{X_1} \cong \PP^3 \\
\cM_{X_1}(x) }
$$
with $\varphi(z) = [E_z]$ and $\pi(z) = F_z$.
We note that there is an isomorphism $\varphi^{-1}(E) \cong \PP^1$ and that  
$\pi(\varphi^{-1}(E)) \subset \cM_{X_1}(x) \subset \PP^5$ is a conic 
for any stable $E \in \cM^s_{X_1}$.
On $X_1 \times Z$ there exists a ``universal'' bundle, which we denote by  
$\cV$ (see \cite{Mu} (3.8)). It has the 
property 
$$\cV_{|X_1 \times \{z\}} \cong E_z, \qquad  \forall z \in Z.$$ 

Let  $\cL$ denote a Poincar\'e bundle on $X_1 \times JX_1$.
By abuse of notation we also denote by $\cV$ and $\cL$ their pull-backs 
to $X_1 \times JX_1 \times Z$. We denote by $\pi_{X_1}$ and $q$ the canonical projections
$$
X_1 \stackrel{\pi_{X_1}}{\longleftarrow} X_1 \times JX_1 \times Z \stackrel{q}{\lra} 
JX_1 \times Z.
$$

We consider the map $m$ given by tensor product
$$m: JX_1 \times \cM_{X_1} \lra \cN_{X_1}, \qquad (L,E) \lms L \otimes E.$$
Note that the restriction of $m$ to the stable locus $m^s : 
JX_1 \times \cM^s_{X_1} \lra \cN^s_{X_1}$ is an \'etale map of degree $16$. We 
denote by $\psi$ the composite map
$$ \psi: JX_1 \times Z \stackrel{\mathrm{id}_{JX_1} \times \varphi}{\lra}
JX_1 \times \cM_{X_1} \stackrel{m}{\lra} \cN_{X_1}, \qquad \psi(L,z) = L
\otimes E_z$$

Let $D \in |\omega_{X_1}|$ be a smooth canonical divisor on 
$X_1$. We introduce the following sheaves over $JX_1 \times Z$
$$\cF_1 = q_*(\cL^* \otimes \cV^* \otimes \pi_{X_1}^*(F_*(\gt^{-1}) \otimes \omega_{X_1})) \qquad
\text{and} \qquad \cF_0 = \oplus_{y \in D} \left( \cL^* \otimes \cV^*_{|\{y\} \times
JX_1 \times Z} \right) \otimes k^{\oplus p}.$$

The next proposition is an even degree analogue of \cite{LN}
Theorem 3.1 . 

\bigskip
 
{\bf Proposition 5.1} {\it
\begin{itemize}
\item[(a)] The sheaves $\cF_0$ and $\cF_1$ are locally free of rank $4p$ and
$4p-4$ respectively and there is an exact sequence
$$ 0 \lra \cF_1 \stackrel{\gamma}{\lra} \cF_0 \lra 
R^1q_*( \cL^* \otimes \cV^* \otimes \pi_{X_1}^*(F_*(\gt^{-1}))) \lra 0.$$
Let $\tilde{\cD} \subset JX_1 \times Z$ denote the subscheme defined by the
$4$-th Fitting ideal of the sheaf $R^1q_*( \cL^* \otimes \cV^* \otimes 
\pi_{X_1}^*(F_*(\gt^{-1})))$. We have set-theoretically
$$
\mathrm{supp} \  \tilde{\cD} = \{ (L,z) \in JX_1 \times Z \ | \ \dim \mathrm{Hom}
(L \otimes E_z, F_*(\gt^{-1})) = 1 \},
$$
and $\dim \tilde{\cD} = 1$.
\item[(b)] Let $\delta$ denote the $l$-adic ($l \not= p$) cohomology class of 
$\tilde{\cD}$ in $JX_1 \times Z$. Then
$$ \delta = c_5(\cF_0 - \cF_1) \in H^{10}(JX_1 \times Z, \ZZ_l).$$

\end{itemize}}

\bigskip

{\it Proof}: We consider the canonical exact sequence over $X_1 \times JX_1 \times Z$ associated
to the effective divisor $\pi_{X_1}^* D$
$$
0 \ra \cL^* \otimes \cV^* \otimes \pi_{X_1}^*F_*(\gt^{-1}) \stackrel{\otimes D}{\lra} 
\cL^* \otimes \cV^* \otimes \pi_{X_1}^*(F_*(\gt^{-1}) 
\otimes \omega_{X_1}) \ra \cL^* \otimes \cV^*_{| \pi_{X_1}^*D} \otimes k^{\oplus p} \ra 0.
$$ 
By Proposition 1.2 the rank-$p$ vector bundle $F_*(\gt^{-1})$ is stable and since 
$$1-\frac{2}{p} = \mu(F_*(\gt^{-1})) > \mu(L \otimes E) = 0 \qquad \forall (L,E) 
\in JX_1 \times \cM_{X_1},$$
we obtain
$$ 
\dim H^1(L^* \otimes E^* \otimes F_*(\gt^{-1}) \otimes \omega_{X_1}) 
= \dim \Hom(F_*(\gt^{-1}), L \otimes E)
= 0. 
$$
This implies
$$
R^1q_{*}(\cL^* \otimes \cV^* \otimes \pi_{X_1}^*(F_*(\gt^{-1}) \otimes \omega_{X_1})) = 0.
$$
By the base change theorems the sheaf $\cF_1$ is locally free. 
Taking direct images by $q$ (note that
$q_{*}(\cL^* \otimes \cV^* \otimes \pi_{X_1}^*F_*(\gt^{-1})) =0$ because it is a torsion sheaf),
we obtain the exact sequence
$$ 0 \lra \cF_1 \stackrel{\gamma}{\lra} \cF_0 \lra 
R^1q_*( \cL^* \otimes \cV^* \otimes \pi_{X_1}^*(F_*(\gt^{-1}))) \lra 0.$$
with $\cF_1$ and $\cF_0$ as in the statement of the proposition. Note that by Riemann-Roch we have
$$
\rk\,\ \cF_1 = 4p - 4 \qquad \mbox{and} \qquad \rk \,\ \cF_0 = 4p.
$$
It follows from the proof of Proposition 3.1 (a) that any nonzero homomorphism 
$L \otimes E \lra F_*(\gt^{-1})$ is injective. Moreover by Proposition 3.1 (b) (iii) 
for any subbundle $L \otimes E \subset F_*(\gt^{-1})$ we
have $\dim \Hom(L \otimes E, F_*(\gt^{-1})) = 1$, or equivalently $\dim
H^1( L^* \otimes E^* \otimes F_*(\gt^{-1})) = 5$. Using the base change theorems we obtain the
following series of equivalences
\begin{align*}
(L,z)  \in \mathrm{supp}  \ \tilde{\cD} & \ \iff \ \rk \ \gamma_{(L,z)} < 4p = \rk \ \cF_0  \\
                              & \ \iff \ \dim H^1( L^* \otimes E^* \otimes F_*(\gt^{-1})) \geq 5 \\
                              & \ \iff \ \dim \Hom(L \otimes E, F_*(\gt^{-1})) \geq 1 \\
                              & \ \iff \ \dim \Hom(L \otimes E, F_*(\gt^{-1})) = 1.
\end{align*}
Finally we clearly have the equality $\mathrm{supp} \  \psi(\tilde{\cD}) = \mathrm{supp} \
\cQ$. Since $\dim \cQ = 0$ and since the fibers of the morphism $\varphi$ over stable vector
bundles are projective lines, we deduce that $\dim \tilde{\cD} = 1$. This proves part (a).

\bigskip

Part (b) follows from Porteous' formula, which says that the fundamental class $\delta
\in H^{10}(JX_1 \times Z, \ZZ_l)$ of the determinantal subscheme $\tilde{\cD}$ 
is given (with the notation of [ACGH], p.86) by
\begin{align*}
\delta  &= \Delta_{4p-(4p-5),4p-4-(4p-5)}(c_t(\cF_0 -\cF_1))\\
                                         &= \Delta_{5,1}(c_t(\cF_0 - \cF_1))\\
                                         &= c_5(\cF_0 - \cF_1). 
\end{align*}   \qed                                    
\bigskip

Let $M$ be a sheaf over a $k$-scheme $S$. We denote by 
$$\mathrm{Fitt}_n[M] \subset \cO_S$$
the $n$-th Fitting ideal sheaf of $M$.

\bigskip

We now define the $0$-dimensional subscheme $\cD \subset \cN_{X_1}^s$, which is supported 
on $\mathrm{supp} \ \cQ$. Consider a
bundle $E \in \cN_{X_1}^s$ with $\dim \mathrm{Hom}(E,F_*(\gt^{-1})) \geq 1$ or
equivalently $\dim H^1(E^* \otimes F_*(\gt^{-1})) \geq 5$. The GIT-construction of the
moduli space  $\cN_{X_1}^s$ realizes  $\cN_{X_1}^s$ as a quotient of an open subset $\cU$ of a Quot-scheme by the group $\PP GL(N)$ for some $N$. It can be shown (see e.g. \cite{La2} section 
3) that $\cU$ is a principal $\PP GL(N)$-bundle for the \'etale topology over $\cN_{X_1}^s$.
Hence there exists an \'etale neighbourhood $\tau : \overline{U} \ra U$ of $E$ over which
the $\PP GL(N)$-bundle is trivial, i.e., admits a section. The universal bundle over the
Quot-scheme restricts to a bundle $\cE$ over $X_1 \times \overline{U}$. Choose a
point $\overline{E} \in \overline{U}$ over $E$. Since $\tau$ is \'etale, it
induces an isomorphism of the local rings $\cO_{\overline{U},\overline{E}}$ and
$\cO_{U,E}$. We simply denote this ring by $\cO_E$. Consider the scheme structure at $\overline{E}$ defined by the ideal $\mathrm{Fitt}_4[R^1 \pi_{\overline{U} *} (\cE^* \otimes \pi_{X_1}^*
F_*(\gt^{-1}))]$. This also defines a scheme structure at $E \in U$, which does
not depend on the choice of the \'etale neighbourhood. Note that there is a
``universal'' bundle $\cE$ over $X_1 \times \mathrm{Spec}(\cO_E)$.

\bigskip

{\bf Lemma 5.2} \ {\it
There is a scheme-theoretical equality}
$$ \tilde{\cD} = \psi^{-1} \cD.$$

\bigskip

{\it Proof:}  We consider a point $E \in \mathrm{supp} \cD = \mathrm{supp} \cQ$ and denote 
by $Z_E$  the fibre $\psi^{-1} (\mathrm{Spec}(\cO_E))$ and
by $\psi_E: Z_E \ra  \mathrm{Spec}(\cO_E)$ the 
morphism obtained from $\psi$ after taking the base change $\mathrm{Spec}(\cO_E) \ra
\cN_{X_1}^s$. The lemma now follows because the formation of the Fitting ideal and
taking the higher direct image $R^1 \pi_{\mathrm{\Spec}(\cO_E) *}$ commutes 
with the base change $\psi_E$ (see \cite{E} Corollary 20.5 and \cite{Ha} 
Proposition 12.5), i.e.
$$ \psi_E^{-1} \left[ \mathrm{Fitt}_4 (R^1\pi_{\mathrm{\Spec}(\cO_E)*}(\cE^*
\otimes \pi_{X_1}^* F_*(\gt^{-1}))  \right] \cdot \cO_{Z_E} = 
\mathrm{Fitt}_4 (R^1\pi_{Z_E*}((\mathrm{id}_{X_1} \times \psi_E)^* \cE^* 
\otimes \pi_{X_1}^* F_*(\gt^{-1})), $$
and $(\mathrm{id}_{X_1} \times \psi_E)^* \cE \sim \cL \otimes \cV_{|X_1 \times Z_E}$.
 
\qed 

\bigskip

{\bf Lemma 5.3} \ {\it
The subscheme $\cD \subset \cN^s_{X_1}$  corepresents the functor
which associates to any $k$-scheme $S$ the set}
\begin{align*}
\underline{\cD}(S) = \{  & \cE \ \text{locally free sheaf over} \ X_1 \times S \ 
\text{of rank} \ 2 \ | \ \deg \cE_{|X_1 \times \{ s \}} = 0 \ \forall s \in S, \\  
& \mathrm{Fitt}_4 [R^1\pi_{S*}( \cE^* \otimes \pi_{X_1}^*(F_*(\gt^{-1})))] = 0 \} /
\sim 
\end{align*}

\bigskip

{\it Proof:}  This is an immediate consequence of the definition of $\cD$ and
the fact that $\cN_{X_1}^s$ universally corepresents the functor
$\underline{\cN}_{X_1}^s$.
\qed

\bigskip

{\bf Lemma 5.4} \ {\it Let $S$ be a $k$-scheme and $\cE$ a sheaf over $X_1 \times
S$ with $\langle \cE \rangle \in \cN_{X_1}^s (S)$. We suppose that the set-theoretical
image of the classifying morphism of  $\cE$ 
$$ \Phi_\cE : S \lra \cN_{X_1}^s, \qquad s \lms \cE_{| X_1 \times \{ s\} } $$
is a point. Then there exists an Artinian ring $A$, a morphism $\varphi : S \lra 
S_0 := \mathrm{Spec}(A)$  and a locally free sheaf $\cE_0$ over  $X_1 \times S_0$ 
such that
\begin{enumerate}
\item $\cE \sim (\mathrm{id}_{X_1} \times \varphi)^* \cE_0$
\item the natural map $\cO_{S_0} \lra \varphi_* \cO_S$ is injective.
\end{enumerate}
}

\bigskip

{\it Proof:}  Since the set-theoretical support of $\mathrm{im} \ \Phi_\cE$ is a 
point, there exists an Artinian ring $A$ such that $\Phi_\cE$ factorizes  through
the inclusion $\mathrm{Spec}(A)  \hookrightarrow \cN_{X_1}^s$. By the 
argument, which we already used in the definition of $\cD$, there exists
a universal bundle $\cE_0$ over $X_1 \times \mathrm{Spec}(A)$. So we have shown 
property (1). As for (2), we consider the ideal $I \subset A$ defined by
$\tilde{I} = \ker ( \cO_{\mathrm{Spec}(A)} \ra \varphi_* \cO_S)$, where $\tilde{I}$
denotes the associated $\cO_{\mathrm{Spec} (A)}$-module. If $I \not= 0$, we 
replace $A$ by $A/I$ and we are done.
\qed 

\bigskip

{\bf Proposition 5.5} \ {\it
There is a scheme-theoretical equality}
$$ \cD = \cQ.$$

\bigskip

{\it Proof:} We note that
$\underline{\cD}(S)$ and $\underline{\cQ}(S)$ are subsets of $\underline{\cN}^s_{X_1}(S)$
(the injectivity of the map $\underline{\cQ}(S) 
\ra \underline{\cN}^s_{X_1}(S)$ is proved similarly as in the proof of Proposition 4.5).
Since $\cD$ and $\cQ$ corepresent the two functors $\underline{\cD}$ and
$\underline{\cQ}$, it will be enough to show that the set $\underline{\cD}(S)$ 
coincides with $\underline{\cQ}(S)$ for any $k$-scheme. 
\bigskip

We first show that $\underline{\cD}(S) \subset \underline{\cQ}(S)$. Consider a sheaf
$\cE$ with $\langle \cE \rangle \in \underline{\cD}(S)$. For simplicity we denote the
sheaf $\cE^* \otimes \pi_{X_1}^* (F_*(\gt^{-1}))$ by $\cH$. By \cite{Ha} Theorem 12.11 
there is an isomorphism
$$ R^1 \pi_{S*} \cH \otimes k(s) \cong H^1(X_1 \times {s}, \cH_{|X_1 \times \{ s\}} )
\qquad \forall s \in S.$$
Since we have assumed $\mathrm{Fitt}_4 [ R^1\pi_{S*} \cH ] = 0$, we obtain 
$\dim H^1(X_1 \times \{ s \}, \cH_{|X_1 \times \{ s \}})  \geq 5$, or
equivalently $\dim H^0(X_1 \times \{ s \}, \cH_{|X_1 \times \{ s \}})  \geq 1$, i.e.,
the vector bundle  $\cE_{|X_1 \times \{ s \}}$ is a subsheaf, hence subbundle, of
$F_*(\gt^{-1})$. This implies that the set-theoretical image of the 
classifying map $\Phi_\cE$ is contained in $\mathrm{supp} \cQ$. Taking connected
components of $S$, we can assume that the image of $\Phi_\cE$ is a point. Therefore
we can apply Lemma 5.4: there exists a locally free sheaf $\cE_0$ over $X_1 \times S_0$
such that $\cE \sim (\mathrm{id}_{X_1} \times \varphi)^* \cE_0 $. For simplicity 
we write $\cH_0 = \cE_0^* \otimes \pi_{X_1}^* (F_*(\gt^{-1}))$. In particular
$\cH = (\mathrm{id}_{X_1} \times  \varphi)^* \cH_0$. Since the projection $\pi_{S_0} : X_1 \times S_0 \ra
S_0$ is of relative dimension $1$, taking the higher direct image $R^1\pi_{S_0*}$
commutes with the (not necessarily flat) base change $\varphi : S \ra S_0$ (\cite{Ha}
Proposition 12.5),
i.e., there is an isomorphism
$$ \varphi^* R^1\pi_{S_0*} \cH_0  \cong R^1\pi_{S*} \cH.$$
Since the formation of Fitting ideals also commutes with any base change (see
\cite{E} Corollary 20.5), we obtain
$$ \mathrm{Fitt}_4[R^1\pi_{S*} \cH] = \mathrm{Fitt}_4 [R^1\pi_{S_0*} \cH_0] \cdot
\cO_S.$$
Since $\mathrm{Fitt}_4[R^1\pi_{S*} \cH] = 0$ and $\cO_{S_0} \ra \varphi_* \cO_S$ is
injective, we deduce that $\mathrm{Fitt}_4[R^1\pi_{S_0*} \cH_0] = 0$. Since 
by Proposition 3.1 (b) (iii) $\dim  R^1 \pi_{S_0*} \cH_0 \otimes k(s_0) = 5$ for the closed point
$s_0 \in S_0$, we have 
$\mathrm{Fitt}_5 [ R^1\pi_{S_0*} \cH_0 ] = \cO_{S_0}$.
We deduce by \cite{E}
Proposition 20.8 that the sheaf $R^1\pi_{S_0*} \cH_0$ is a free $A$-module of rank $5$. 
By \cite{Ha} Theorem 12.11 (b) we deduce that there is an isomorphism
$$ \pi_{S_0*} \cH_0 \otimes k(s_0) \cong H^0(X_1 \times {s_0}, \cH_{|X_1 \times \{ s_0\}} )
$$
Again by Proposition 3.1 (b) (iii) we obtain $\dim \pi_{S_0*} \cH_0 \otimes k(s_0) = 1$.
In particular the $\cO_{S_0}$-module $\pi_{S_0*} \cH_0$ is not zero and therefore there
exists a nonzero global section $i \in H^0(S_0,\pi_{S_0*} \cH_0) = H^0(X_1 \times S_0,
\cE_0^* \otimes \pi^*_{X_1} F_*(\gt^{-1}))$. We pull-back $i$ under the map
$\mathrm{id}_{X_1} \times \varphi$ and we obtain a nonzero section
$$ j = (\mathrm{id}_{X_1} \times \varphi)^* i \in H^0(X_1 \times S,
\cE^* \otimes \pi^*_{X_1} F_*(\gt^{-1})).$$
Now we apply Lemma 4.3 and we continue as in the proof of Proposition 4.6. This shows  
that $\langle \cE \rangle \in \underline{\cQ}(S)$.

\bigskip

We now show that $\underline{\cQ}(S) \subset \underline{\cD}(S)$. Consider a sheaf
$\cE \in \underline{\cQ}(S)$. The nonzero global section $j \in 
H^0(X_1 \times S, \cH) = H^0(S, \pi_{S*} \cH)$ determines by evaluation
at a point $s \in S$ an element $\alpha \in \pi_{S*} \cH \otimes k(s)$. The 
image of $\alpha$ under the natural map 
$$ \varphi^0(s) : \pi_{S*}\cH \otimes k(s) \lra H^0(X_1 \times \{ s \}, 
\cH_{| X_1 \times \{ s \} } ) $$
coincides with $j_{| X_1 \times \{ s \}}$ which is nonzero. Moreover since
$\dim H^0(X_1 \times \{ s \}, \cH_{| X_1 \times \{ s \} }) = 1$, we obtain that
$\varphi^0(s)$ is surjective. Hence by \cite{Ha} Theorem 12.11 the sheaf
$R^1 \pi_{S*} \cH$ is locally free of rank $5$. Again by \cite{E} Proposition
20.8 this is equivalent to $\mathrm{Fitt}_4[R^1 \pi_{S*} \cH] = 0$ and
$\mathrm{Fitt}_5[R^1 \pi_{S*} \cH] = \cO_S$ and we are done.  \qed

\bigskip

\begin{center}
{\bf \S 6 Chern class computations.}
\end{center}
\vspace*{0.5cm}

In this section we will compute the length of the determinantal subscheme
$\cD \subset \cN_{X_1}$ by evaluating the Chern 
class $c_5(\cF_0 - \cF_1)$ --- see Proposition 5.1 (b). 

\bigskip

Let $l$ be a prime number different from $p$. We have to recall some properties of the 
cohomology ring $H^*(X_1 \times JX_1 \times Z, \ZZ_l)$ (see also \cite{LN}).    
In the sequel we identify all classes of $H^*(X_1,\ZZ_l), \,\,\, H^*(JX_1,\ZZ_l)$ etc. 
with their preimages in $H^*(X_1 \times JX_1 \times Z, \ZZ_l)$ under 
the natural pull-back maps. 

\bigskip

Let $\Theta \in H^2(JX_1,\ZZ_l)$ denote the class of the theta divisor in $JX_1$. 
Let $f$ denote a positive generator of 
$H^2(X_1,\ZZ_l)$. The cup product  
$H^1(X_1,\ZZ_l) \times H^1(X_1,\ZZ_l) \ra H^2(X_1,\ZZ_l) \simeq \ZZ_l$ 
gives a symplectic structure on 
$H^1(X_1,\ZZ_l)$. Choose a symplectic basis $e_1,e_2,e_3,e_4$ of $H^1(X_1,\ZZ_l)$ such 
that $e_1e_3 = e_2e_4 = -f$ and all other products $e_ie_j = 0$.
We can then normalize the Poincar\'e bundle $\cL$ on $X_1 \times JX_1$ so that
\begin{equation} \label{e1}
c(\cL) = 1 + \xi_1   
\end{equation}
where $\xi_1 \in H^1(X_1,\ZZ_l) \otimes H^1(JX_1,\ZZ_l) \subset H^2(X_1 \times JX_1,\ZZ_l)$ 
can be written as 
$$
\xi_1 = \sum_{i=1}^4 e_i \otimes \varphi_i
$$
with $\varphi_i \in H^1(JX_1,\ZZ_l)$. Moreover, we have by the same reasoning, applying
[ACGH] p.335 and p.21 
\begin{equation} \label{e2}
\xi_1^2 = -2\Theta f \quad \mbox{and} \quad  \Theta^2[JX_1] =2.   
\end{equation}
Since the variety $\cM_{X_1}(x)$ is a smooth intersection of 2 quadrics in $\PP^5$, one can work out that the $l$-adic cohomology groups $H^i(\cM_{X_1}(x),\ZZ_l)$ for $i=0,\ldots, 6$ are
(see e.g. \cite{Re} p. 0.19)
$$
\ZZ_l, \,\,\,0, \,\,\,\ZZ_l,\,\,\, \ZZ^4_l,\,\,\, \ZZ_l, \,\,\,0, \,\,\,\ZZ_l.
$$
In particular $H^2(\cM_{X_1}(x), \ZZ_l)$ is free of rank $1$ and, if
$\alpha$ denotes a positive generator of it, then
\begin{equation} \label{eq1}
\alpha^3[\cM_{X_1}(x)] = 4.   
\end{equation}
According to [N2] p. 338 and applying 
reduction mod p and a comparison theorem, the Chern classes of the universal bundle $\cU$ are of the form
\begin{equation} \label{e3}
c_1(\cU) = \alpha + f \qquad \mbox{and} \qquad c_2(\cU) = \chi + \xi_2 +\alpha f 
\end{equation}
with $\chi \in H^4(\cM_{X_1}(x),\ZZ_l)$ and $\xi_2 \in H^1(X_1,\ZZ_l) \otimes 
H^3(\cM_{X_1}(x),\ZZ_l)$.
As in \cite{N} and \cite{KN} we write 
\begin{equation} \label{eqn}
\beta = \alpha^2 - 4\chi \quad \mbox{and} \quad \xi_2^2 = \gamma f \quad \mbox{with} \quad \gamma \in H^6(\cM_{X_1}(x),\ZZ_l).
\end{equation}   
Then the relations of \cite{KN} give
$$
\alpha^2 + \beta = 0 \quad \mbox{and} \quad \alpha^3 + 5 \alpha \beta + 4 \gamma = 0.
$$
Hence $\beta = - \alpha^2, \; \gamma = \alpha^3$. Together with (\ref{e3}) and (\ref{eqn}) this gives 
\begin{equation} \label{eqn1}
c_2(\cU) = \frac{\alpha^2}{2} + \xi_2 + \alpha f \quad \mbox{and} \quad \xi_2^2 = \alpha^3f
\end{equation}
Define $\Lambda \in H^1(JX_1,\ZZ_l) \otimes H^3(\cM_{X_1}(x),\ZZ_l)$ by
\begin{equation} \label{e5}
\xi_1 \xi_2 = \Lambda f.
\end{equation}
Then we have for dimensional reasons and noting that $H^5(\cM_{X_1}(x),\ZZ_l) = 0$, that the following classes are all zero:
\begin{equation} \label{e6}
f^2, \,\, \xi_1^3, \,\, \alpha^4, \,\, \xi_1 f, \,\, \xi_2 f,\,\, \alpha \xi_2, \,\, \alpha \Lambda,\,\, \Theta^2\Lambda,\,\, \Theta^3.   
\end{equation}
Finally, $Z$ is the $\PP^1$-bundle associated to the vector bundle $\cU_x$ on 
$\cM_{X_1}(x)$. Let $H \in H^2(Z,\ZZ_l)$ denote the
first Chern class of the tautological line bundle on $Z$. We have, using the definition of the Chern classes $c_i(\cU)$
 and (\ref{eq1}), 
\begin{equation} \label{e7}
H^2 = \alpha H - \frac{\alpha^2}{2}, \qquad H^4 = 0, \qquad \alpha^3 H[Z] = 4  
\end{equation}
and we get for the ``universal'' bundle $\cV$,
\begin{equation} \label{e8}
c_1(\cV) = \alpha \quad \mbox{and} \quad c_2(\cV) = \frac{\alpha^2}{2} + \xi_2 + Hf.  
\end{equation}

\bigskip

{\bf Lemma 6.1} {\it 
\begin{enumerate}
\item[(a)] The cohomology class $\alpha \cdot c_5(\cF_0 - \cF_1) 
\in H^{12}(JX_1 \times Z, \ZZ_l)$ is a multiple of  the class
$\alpha^3H\Theta^2$. 
\item[(b)] The pull-back under the map $\varphi: Z \lra \cM_{X_1} \cong \PP^3$ of the class of a point is the class $H^3 = \frac{\alpha^2}{2}H - \frac{\alpha^3}{2}$.
\end{enumerate}}

\bigskip

{\it Proof}: For part (a)  it is enough to note that all other relevant cohomology classes vanish, 
since $\alpha^4 = 0$ and $\alpha\Lambda = 0$. \\
As for part (b), it suffices to show that $c_1(\varphi^* \cO_{\PP^3}(1)) = H$. The line
bundle $\cO_{\PP^3}(1)$ is the inverse of the determinant line bundle \cite{KM} over the
moduli space $\cM_{X_1}$. Since the formation of the determinant line bundle commutes
with any base change (see \cite{KM}), the pull-back  $\varphi^* \cO_{\PP^3}(1)$ is the
inverse of the determinant line bundle associated to the family $\cV \otimes \pi_{X_1}^* N$
for any line bundle $N$ of degree $1$ over $X_1$. Hence the first Chern class of 
$\varphi^* \cO_{\PP^3}(1)$ can be computed by the Grothendieck-Riemann-Roch theorem
applied to the sheaf $\cV \otimes \pi_{X_1}^* N$ over $X_1 \times Z$ and the 
morphism $\pi_Z : X_1 \times Z \ra Z$. We have
\begin{align*}
ch(\cV \otimes \pi_{X_1}^* N) \cdot \pi_{X_1}^* td(X_1) & = \left(2 + \alpha + (-\xi_2 - Hf)
+ \mathrm{h.o.t.} \right) (1+f) (1-f) \\
 & = 2 + \alpha + (-\xi_2 - Hf) + \mathrm{h.o.t.},
\end{align*}
and therefore G-R-R implies that $c_1(\varphi^* \cO_{\PP^3}(1)) = H$ --- note 
that $\pi_{Z*}(\xi_2) = 0$.
 \qed\\

\bigskip

{\bf Proposition 6.2}  \ {\it We have  }
$$
l(\cD) = \frac{1}{24}p^3(p^2-1).
$$

\bigskip

{\it Proof}: Let $\lambda$ denote the length of the subscheme $m^{-1}(\cD)
\subset JX_1 \times \cM_{X_1}$ Since the map $m^s$ is \'etale  of degree $16$, 
we obviously have the relation $\lambda = 16 \cdot l(\cD)$.
According to Lemma 6.1 (b) we have in $H^{10}(JX_{1} \times Z,\ZZ_l)$ 
$$
[(\id \times \varphi)^{-1}(pt)] = H^3\cdot \frac{\Theta^2}{2} = 
\frac{1}{4}\alpha^2H\Theta^2 - \frac{1}{4}\alpha^3\Theta^2,
$$
where $pt$ denotes the class of a point in $JX_1 \times \cM_{X_1}$.
Using Lemma 5.2 we obtain that 
the class $\delta = c_5(\cF_0 - \cF_1) \in H^{10}(JX_{1} \times Z, \ZZ_l)$ equals
$\lambda \cdot (\frac{1}{4}\alpha^2H\Theta^2 - \frac{1}{4}\alpha^3\Theta^2)$. 
Intersecting with $\alpha$ we obtain with Lemma 6.1 (a) and \eqref{e6}
\begin{equation} \label{deflam}
 \alpha \cdot c_5(\cF_0 - \cF_1) = 
\frac{\lambda}{4} \alpha^3H\Theta^2.
\end{equation}
So we have to compute the class $\alpha \cdot c_5(\cF_0 - \cF_1)$. 
By \eqref{e1} and \eqref{e2}, 
$$
ch(\cL) = 1 + \xi_1 - \Theta f
$$
whereas by \eqref{eqn1}, \eqref{e6} and \eqref{e8},
$$
ch(\cV) = 2 + \alpha + (-\xi_2 - Hf) + \frac{1}{12}(-\alpha^3 - 6\alpha Hf) + 
\frac{1}{12}(\alpha^3f - \alpha^2Hf). 
$$
Moreover 
$$
ch(\pi_{X_1}^*(F_*(\gt^{-1}) \otimes \omega_{X_1})) \cdot \pi_{X_1}^*td(X_1) = p + (2p - 2)f.
$$
So using \eqref{eqn1}, \eqref{e5} and \eqref{e6},
\begin{align*}
ch (\cV^* \otimes \cL^* \otimes \,\,\, & \pi_{X_1}^* (F_*(\theta^{-1}) \otimes \omega_{X_1})) \cdot \pi_{X_1}^* td (X_1) 
= 2p +[(4p-4)f - p \alpha - 2p\xi_1]\\
& + \left[p\alpha \xi_1 - 2p\Theta f - (2p-2)\alpha f -p \xi_2 - pHf \right]\\
& + \left[\frac{p}{12}\alpha^3 + \frac{p}{2} \alpha Hf + p \Lambda f + p \alpha \Theta f \right]\\
& + \left[\frac{3p - 2}{12} \alpha^3 f - \frac{p}{12} \alpha^3 \xi_1 - \frac{p}{12} \alpha^2 Hf \right] + 
\left[- \frac{p}{12} \alpha^3 \Theta f \right].\\
\end{align*}
Hence by Grothendieck-Riemann-Roch
for the morphism $q$ we get
\begin{align*}
ch (\cF_1) = 4p-4 &+ [-(2p-2)\alpha -2p\Theta -pH] +\left[ \frac{p}{2} \alpha H + p \Lambda + p \alpha \Theta \right]\\
& +\left[ \frac{3p-2}{12} \alpha^3 - \frac{p}{12} \alpha^2 H \right] + \left[ - \frac{p}{12} \alpha^3 \Theta \right]. \\
\end{align*}
From \eqref{e2} and \eqref{e8} we easily obtain
$$
ch(\cF_0) = 4p - 2p \alpha + \frac{p}{6} \alpha^3.
$$
So
\begin{align*}
ch(\cF_0 - \cF_1) = 4 & + [2p\Theta - 2\alpha +pH] + \left[ - \frac{p}{2} \alpha H - p \Lambda - p \alpha \Theta \right]\\
& + \left[ - \frac{p+1}{12} \alpha^3 + \frac{p}{12}\alpha^2 H \right] + \left[ \frac{p}{12} \alpha^3 \Theta \right].
\end{align*}

%
%Now we are in a position to compute $\lambda$:\\
\noindent
Defining $p_n := n! \cdot ch_n(\cF_0 -\cF_1)$ we have according to Newton's recursive formula ([F] p.56),
$$
c_5(\cF_0 -\cF_1) = \frac{1}{5} \left( p_5 - \frac{5}{6}p_2p_3 - \frac{5}{4}p_1p_4 + \frac{5}{6}p_1^2p_3 + \frac{5}{8}p_1 p_2^2 
- \frac{5}{12}p_1^3p_2 + \frac{1}{24}p_1^5 \right)
$$
with
\begin{align*}
p_1 &= 2p\Theta - 2\alpha + pH\\
p_2 &= -p(\alpha H + 2\Lambda + 2\alpha \Theta)\\
p_3 &= \frac{1}{2}(-(p+1)\alpha^3 + p \alpha^2 H) \\
p_4 &= 2p \alpha^3 \Theta\\
p_5 &= 0.
\end{align*} 
Now an immediate computation using \eqref{e6} and \eqref{e7} gives
$$
\alpha \cdot c_5(\cF_0  - \cF_1) = \frac{p^3(p^2 - 1)}{6} \alpha^3 H \Theta^2.
$$
We conclude from \eqref{deflam} that $\lambda = \frac{2}{3} p^3(p^2 -1)$ and we are 
done. \qed \\

\bigskip

{\bf Remark 6.3}  If $k = \CC$, the number of maximal subbundles of a general vector 
bundle has recently been computed by Y. Holla by using Gromov-Witten invariants \cite{Ho}.
His formula (\cite{Ho} Corollary 4.6) coincides with ours.

\bigskip

\begin{center}
{\bf \S 7 Proof of Theorem 2.}
\end{center}
\bigskip

The proof of Theorem 2 is now straightforward. It suffices to combine
Corollary 4.7, Proposition 5.5 and Proposition 6.2 to obtain the length
$l(\cB)$.

\bigskip

The fact that $\cB$ is a local complete intersection follows from the
isomorphism $\cB_\gt = \cQ_0$ (Proposition 4.6) and Proposition 4.1. \qed

\bigskip

\begin{center}
{\bf \S 8 Questions and Remarks.}
\end{center}
\bigskip

\begin{enumerate}
\item Is the rank-$p$ vector bundle $F_* L$ very stable, i.e. $F_*L$ has no
nilpotent $\omega_{X_1}$-valued endomorphisms, for a general line bundle?

\bigskip

\item Is $F_*(\gt^{-1})$ very stable for a general curve $X$? Note that
very-stability of $F_*(\gt^{-1})$ implies reducedness of $\cB$ (see e.g.
\cite{LN} Lemma 3.3).

\bigskip

\item If $g=2$, we have shown that for a general stable $E \in \cM_{X}$ the
fibre $V^{-1}(E)$ consists of $\frac{1}{3} p (p^2+2)$  stable vector bundles
$E_1 \in \cM_{X_1}$, i.e. bundles $E_1$ such that $F^*E_1 \cong E$ or 
equivalently (via adjunction) $E_1 \subset F_* E$. The Quot-scheme
parametrizing rank-$2$ subbundles of degree $0$ of the rank-$2p$ vector
bundle $F_*E$ has expected dimension $0$, contains the fibre $V^{-1}(E)$,
but it also has a $1$-dimensional component arising from Frobenius-destabilized
bundles.

\bigskip

\item If $p=3$ the base locus $\cB$ consists of $16$ reduced points, which 
correspond to the $16$ nodes of the Kummer surface associated to $JX$ (see 
\cite{LP2} Corollary 6.6). For general $p$, does the configuration of 
points determined by $\cB$ have some geometric significance? 

\end{enumerate}

\end{document}